\documentclass[10pt]{amsart}
\usepackage{amssymb,amsmath,amscd,amsfonts,amsthm,mathrsfs}
\usepackage{verbatim,paralist}
\input xy
\xyoption{all}
\usepackage[all]{xy}
\usepackage{color}

\usepackage[latin1]{inputenc}
\newcommand{\M}{\mathbb{M}}
\newcommand{\N}{\mathbb{N}}

\newcommand{\Q}{\mathbb{Q}}
\newcommand{\R}{\mathbb{R}}
\newcommand{\C}{\mathbb{C}}
\newcommand{\cc}{C}% connected component
\providecommand{\off}{\operatorname{off}}% offspring
\newcommand{\E}{\mathbb E}% expectation
\providecommand{\Fc}{\mathcal F}% sigma algebra
\newcommand{\Ac}{\mathcal{A}}

\newcommand{\Cc}{\mathcal{C}}
\newcommand{\Dc}{\mathcal{D}}

\newcommand{\Hc}{\mathcal{H}}

\newcommand{\Rc}{\mathcal{R}}

\newcommand{\Vc}{\mathcal{V}}
\def\cchi{\raisebox{.45 ex}{$\chi$}}%an upper placed \chi

\providecommand{\father}[1]{\overleftarrow{#1}}

\newcommand{\Hr}{\mathscr{H}}
\newcommand{\Kr}{\mathscr{K}}
\newcommand{\Mr}{\mathscr{M}}
\newcommand{\Nr}{\mathscr{N}}

\newcommand{\dG}{d_G} % degree in G
\newcommand{\tdG}{\tilde d_G} % degree in G
\newcommand{\dV}{\rho_V} % Metric on V
\newcommand{\sphere}[1]{S_{#1}} % Sphere of radius #1 in a graph
\newcommand{\hatsphere}[1]{\hat S_{#1}}%
\RequirePackage{ifthen}
\newcommand{\norm}[2][]{\lVert#2\rVert\ifthenelse{\equal{#1}{}}{}{_{#1}}}%
\newcommand{\bignorm}[2][]{\bigl\|#2\bigr\|\ifthenelse{\equal{#1}{}}{}{_{#1}}}%
\def\build#1_#2^#3{\mathrel{\mathop{\kern 0pt#1}\limits_{#2}^{#3}}}
\newtheorem{theorem}{Theorem}[section]%{\bf}{\it }
\newtheorem{proposition}{Proposition}[section]%{\bf}{\it }
\newtheorem{lemma}{Lemma}[section]%{\bf}{\it }
\newtheorem{remark}{Remark}[section]%{\bf}{\it }
\newtheorem{example}{Example}[section]%{\bf}{\it }
\newtheorem{definition}{Definition}[section]%{\bf}{\it }
\numberwithin{equation}{section}

\newcommand{\slim}{\operatorname*{s-lim}}
\newcommand{\ran}{\mathrm{ran}}

\def \bone{\mathbf{1}}

\def \rmi{{\rm i}}

\def \im{{\rm Im}}
\title[Deficiency indices for discrete Schr\"odinger operators]
      {The problem of deficiency indices for discrete Schr\"odinger operators
	on locally finite graphs}
\begin{document}
%%%%%%%%%%%%%%%%%%%%%%%AUTHOR%%%%%%%%%%%%%%%%%%%%%%%%%%%%%%%%%%%%%%%%%%%
\author{Sylvain Gol\'enia}
\address{Mathematisches Institut der Universit\"at Erlangen-N\"urnberg,
Bismarckstr.\ 1 1/2 \\
91054 Erlangen, Germany
\\
and
Institut de Math\'ematiques de Bordeaux, Universit\'e
Bordeaux $1$, $351,$ cours de la Lib\'eration
\\$33405$ Talence cedex, France
}
\email{sylvain.golenia@u-bordeaux1.fr}
\author{Christoph Schumacher}
\address{Fakult\"at f\"ur Mathematik der Technischen Universit\"at Chemnitz,
Reichenhainer Str.~41\\
09126 Chemnitz, Germany}
\email{christoph.schumacher@mathematik.tu-chemnitz.de}
%\subjclass[2000]{05C63, 05C50, 47B25, 05C80, 05C05}
\keywords{adjacency matrix, discrete Laplacian, deficiency indices,
  locally finite graphs, self-adjointness, spectral graph theory,
  trees, random trees, discrete Schr\"odinger operators, Carleman condition}
\date{Version of \today}
%Abstract%%%%%%%%%%%%%%%%%%%%%%%%%%%%%%%%%%%%%%%%%%%%%%%%%%%%%%%%%%%%%%%
\begin{abstract}
The number of self-adjoint extensions of a symmetric operator acting
on a complex Hilbert space is characterized by its deficiency indices.
Given a locally finite unoriented simple tree, we prove that
the deficiency indices of any discrete Schr\"odinger operator
are either null or infinite.
We also prove that all deterministic discrete Schr\"odinger operators
which act on a random tree are almost surely selfadjoint.
Furthermore, we provide several criteria of essential
self-adjointness. We also address some importance to the case
of the adjacency matrix and conjecture that, given a locally finite
unoriented simple graph, its deficiency indices are either null or
infinite. Besides that, we consider some generalizations of trees and
weighted graphs.
\end{abstract}

\maketitle
\tableofcontents
\section{Introduction}
The spectral theory of adjacency matrices acting on
graphs is useful for the study, among others, of some gelling
polymers, of some electrical networks, and in number theory,
e.g., \cite{CDS, DS, DSV, MO}. In quantum
physics, proving that a symmetric operator is self-adjoint is a
central problem. To characterize all the possible extensions, one
studies the so-called deficiency indices.

We start with some definitions to fix notation for graphs and refer to
\cite{CdV, Chu, MW} for surveys on the matter.
Let $V$ be a countable set.
We equip $V$ with the discrete topology.
Let $E:=V\times V\rightarrow [0,\infty)$ and assume that $E(x,y)=E(y,x)$,
for all $x,y\in V$. We say that $G:=(E,V)$ is an
unoriented weighted graph with \emph{vertices}~$V$ and
\emph{weights}~$E$. In the setting of electrical networks,
the weights correspond to the conductances.
We say that $x,y\in V$ are \emph{neighbors} if $E(x,y)\neq 0$
and denote it by $x\sim y$.
We say there is a \emph{loop} in $x\in V$ if $E(x,x)\neq 0$.
A graph~$G$ is \emph{simple} if it has no loops and $E$ has values in $\{0,1\}$.
The set of \emph{neighbors} of $x\in E$ is denoted by
$\Nr_G(x):=\{y\in E\mid x\sim y\}$.
Given $X\subseteq V$ we write $\Nr_G(X):=\bigcup_{x\in X}\Nr_G(x)$.
The \emph{degree} of $x\in V$ is by definition $\dG(x):=|\Nr_G(x)|$,
the number of neighbors of~$x$.
The graph is of \emph{bounded degree}, if $\sup_{x\in V}\dG(x)$ is finite.
A graph is \emph{locally finite} if $\dG(x)$ is finite for all $x\in
V$. A graph is \emph{connected},
if for all $x,y\in V$, there exists an $x$-$y$-\emph{path}, i.e.,
there is a finite sequence $(x_1,\dotsc,x_n)\in V^{N+1}$ such that
$x_1=x$, $x_{N+1}=y$ and $x_n\sim x_{n+1}$, for all $n\in\{1,\dotsc,N\}$.
In this case, we endow $V$ with the metric $\dV$ defined by
$\dV(x,y):=\inf\{n\in\N\mid\text{there exists an $x$-$y$-path of length~$n$}\}$.
Note that in this paper we use $\N$ for the set of nonpositive
  integers, i.e., $0\in\N$.
In the sequel, all graphs are supposed to be locally finite,
with no loops and unoriented.

We associate to $G$ the complex Hilbert space $\ell^2(V)$.
We denote by $\langle\cdot,\cdot\rangle$ and by $\norm{\cdot}$
the scalar product and the associated norm, respectively.
By abuse of notation, we denote the space simply by $\ell^2(G)$.
The set of complex functions with compact support in~$V$ is denoted by
$\Cc_c(G)$.
One often considers the \emph{Laplacian} defined by
\begin{eqnarray}\label{e:defreg}
  (\Delta_{G,\circ}f)(x):=\sum_{y\sim x}E(x,y)\big(f(x)-f(y)\big),
    \text{ with }f\in\Cc_c(G)
\end{eqnarray}
and the so-called \emph{adjacency matrix}:
\begin{eqnarray}\label{e:def}
  (\Ac_{G,\circ}f)(x):=\sum_{y\sim x}E(x,y)f(y),\text{ with }f\in\Cc_c(G).
\end{eqnarray}
Both of them are symmetric and thus closable.
We denote the closures by $\Delta_{G}$ and $\Ac_{G}$,
their domains by $\Dc(\Delta_G)$ and $\Dc(\Ac_G)$,
and their adjoints by $(\Delta_G)^*$ and $(\Ac_G)^*$, respectively.
In \cite{Woj}, see also \cite{Jor08}, it is shown that the operator
$\Delta_G$ is essentially self-adjoint on $\Cc_c(G)$, when the graph is simple.
In particular, one has that $\Delta_G=(\Delta_G)^*$.
In contrast, even in the case of a locally finite tree~$G$,
$\Ac_G$ may have many self-adjoint extensions,
see \cite{MO, Mu, Gol} and Proposition \ref{p:tree} for concrete examples.
We mention also the work \cite{Aom}, where a characterization
in terms of limit point -- limit circle is given.

In this note, we are also interested in the
\emph{discrete Schr\"odinger operators}
$\Ac_G+\Vc$ and $\Delta_{G}+\Vc$ with \emph{potential} $\Vc:=V\to\R$,
where $\Vc$ also denotes the operator of multiplication with the function $\Vc$.
The operators are defined as the closures of $\Ac_{G,\circ}+\Vc$
and of $\Delta_{G,\circ}+\Vc$ on $\Cc_c(G)$, respectively.
Note that $\Delta_G$, up to sign, is in fact a discrete Schr\"odinger
operator formed with the help of $\Ac_G$:
\begin{eqnarray}\label{e:relation}
  \Delta_G=\Vc-\Ac_G,\text{ where }\Vc(x):=\sum_{y\sim x}E(x,y).
\end{eqnarray}
In the sequel, we investigate the number of possible self-adjoint
extensions of discrete Schr\"odinger operators by computing their
deficiency indices.  Given a closed and densely defined symmetric operator~$T$
acting on a complex Hilbert space, the deficiency indices of~$T$
are defined by $\eta_\pm(T):=\dim\ker(T^* \mp \rmi)\in\N\cup\{+\infty\}$.
We recall some well-known facts. The operator~$T$ possesses a self-adjoint
extension if and only if $\eta_+(T)=\eta_-(T)$.  If this is the case,
we denote the common value by $\eta(T)$.
$T$ is self-adjoint if and only if $\eta(T)=0$.
Moreover, if $\eta(T)$ is finite,
the self-adjoint extensions can be explicitly parametrized by the unitary
group $U(n)$ in dimension $n=\eta(T)$.  Using the Krein formula,
it follows that the absolutely continuous spectrum of all self-adjoint
extensions is the same.

Since the operator $\Ac_G+\Vc$ commutes with the complex conjugation,
its deficiency indices are equal, e.g., \cite[Theorem X.3]{RS}.
We denote by $\eta(G)$ the common value, when $\Vc=0$.
This means that $\Ac_G+\Vc$ possesses a self-adjoint extension.
Remark that $\eta(\Ac_G+\Vc)=0$ (resp.\ $\eta(\Delta_G+\Vc)=0$)
if and only if $\Ac_G+\Vc$ (resp.\ $\Delta_G+\Vc$)
is essentially self-adjoint on $\Cc_c(G)$.
We give the following criteria for essential self-adjointness:
\begin{proposition}\label{p:essSA}
  Let $G=(E,V)$ be a locally finite graph and $\Vc:V\to\R$ be
  a potential. Then, the following assertions hold true:
  \begin{enumerate}
  \item Provided that $\Vc$ is bounded from below,
    $\Delta_G+\Vc$ is essentially self-adjoint on $\Cc_c(G)$.
  \item Let $x_0\in V$, set
    $b_i:=\sup\{\sum_{x,y}E(x,y)\mid
      \dV(x_0,x)=i\text{ and }\dV(x_0,y)=i+1\}$,
    and take $\Vc:V\to\R$.  If $\sum_{i\in\N}1/b_i=+\infty$,
    then $\Ac_G+\Vc$ and $\Delta_G+\Vc$
    is essentially self-adjoint on $\Cc_c(G)$.
  \item Suppose that $\sup_x\max_{y\sim x}|\dG(x)-\dG(y)|<\infty$,
    $E$ is bounded, and $\sup_{x\in V}|\Vc(x)/\dG(x)|<\infty$,
    then $\Ac_G+\Vc$ is essentially self-adjoint on $\Cc_c(G)$.
  \item Suppose that $\dG$ is bounded,
    $\sup_x\max_{y\sim x}|E(x)-E(y)|<\infty$,
    where $E(x):=\max_{y\sim x}E(x,y)$,
    and that $\sup_{x\in V}|\Vc(x)/E(x)|<\infty$,
    then $\Ac_G+\Vc$ is essentially self-adjoint on $\Cc_c(G)$.
  \item Suppose there is a compact set $K\subset V$, such that
    $\sum_{y\sim x}E^2(x,y)\dG(y)\le\Vc^2(x)$ for all $x\notin K$.
    Then $\Ac_G+\Vc$ is essentially self-adjoint on $\Cc_c(G)$.
  \item Suppose there is a compact set $K\subset V$, such that
    $\sum_{y\sim x}E^2(x,y) \big(1+\dG(y)\big)\leq \Vc^2(x)$ for all
    $x\notin K$, then $\Delta_G+\Vc$ is essentially self-adjoint on
    $\Cc_c(G)$.
  \end{enumerate}
\end{proposition}
We prove the result in Section \ref{s:essSA}. The first point is the
discrete version of the fact that given a
non-negative potential $\Vc\in L^2_{\rm loc}(\R^n)$, one has that
$-\Delta_{\R^n}+\Vc$ is essentially self-adjoint on $\Cc_c(\R^n)$,
e.g., \cite[Theorem X.28]{RS}. It is essentially a repetition of
\cite[Theorem 1.3.1]{Woj}. The second point is a Carleman-type
condition, see for instance \cite[Page 504]{Be} for the case of
Jacobi matrices. We stress that this result holds true without any
hypothesis of size or of sign on the potential part. In particular,
the Schr\"odinger operators could be unbounded from below and from
above, see \cite{Gol} for instance. Unlike in \cite{Be}, we rely on an
commutator approach, see \cite{Wo1, Wo2} for similar techniques.
The points (3) and (4) follow by application of the Nelson commutator
Theorem. The two last ones are an application of W\"ust's Theorem by
considering $\Ac$ and $\Delta$ as perturbation of the potential. We
mention the works of \cite{CTT, KL, Ma} on related questions.

Concentrate a moment on the case of the adjacency matrix for simple
graphs. Keep in mind, it is no gentle perturbation of the Laplacian,
see Proposition \ref{p:noneq}. In \cite{MO, Mu}, adjacency matrices
for simple trees with positive deficiency indices are constructed.
In fact, it follows from the proof that the deficiency indices
are infinite in both references.
We recall that a \emph{tree} is a connected graph $G=(E,V)$
such that for each edge $e\in V\times V$ with $E(e)\ne0$
the graph $(\tilde E, V)$, with $\tilde E := E \times 1_{\{e\}^c}$,
i.e., with $e$ removed, is disconnected.
As a general result, a special case of
Theorem \ref{t:main0} gives that, given a locally finite simple tree~$G$,
one has
\begin{eqnarray}\label{e:problem}
  \eta(G)\in\{0,+\infty\}.
\end{eqnarray}
This is a new result to our knowledge, although
the literature on trees is extensive. We believe that,
given a simple graph $G=(E,V)$,
or more generally, a graph with bounded weights,
\eqref{e:problem} should be true. In Remark \ref{r:tensor},
we explain that it is enough to prove \eqref{e:problem}
for simple bi-partite graphs. We recall that a \emph{bi-partite} graph
is a graph so that its vertex set can be partitioned into two subsets
in such a way that no two points in the same subset are
neighbors. Trees are bi-partite for instance. We stress that this
conjecture is false if one takes unbounded weights, see for instance
counter-examples of adjacency matrices given by Jacobi matrices in
\cite[Remark 2.1]{Gol} and also in \cite{MO}.

We now point out that the self-adjointness of the adjacency matrix,
acting on a simple locally finite tree~$G$,
is linked with the growth of the offspring, i.e., of the number of
sons. (We refer to Section~\ref{s:def} for precise definitions
concerning trees.) When the latter grows up to linearly,
Proposition~\ref{p:essSA} gives that $\eta(G)=0$. On the other hand,
if the growth is ``exponential", Proposition~\ref{p:exp} assures that
$\eta(G)=\infty$. In Section~\ref{s:birth}, using invariant spaces,
we prove the following sharp result:
\begin{proposition}\label{p:tree}
  Let $\alpha>0$ and $G$ be a tree with offspring
  $\lfloor{n^\alpha}\rfloor$ per individual at generation~$n$.
  Then, one obtains:
  \begin{eqnarray*}
    \eta(G)=
    \begin{cases}
      0,	&\text{ if }\alpha\le2,\\
      +\infty,	&\text{ if }\alpha>2.
    \end{cases}
  \end{eqnarray*}
\end{proposition}

We come back to the general question for Schr\"odinger operators and
give our main result in the context of trees. We prove it in Section
\ref{s:possible} and generalize it in Theorem \ref{t:main} to
a family of graphs obtained recursively.
\begin{theorem}\label{t:main0}
Let $G=(E,V)$ be a locally finite weighted tree, where $E$ is
bounded, and let $\Vc:V\rightarrow \R$ be a potential. Then one has:
\begin{eqnarray}\label{e:schroe}
\eta(\Ac_G+\Vc)\in \{0, +\infty\} \text{ and } \eta(\Delta_G+\Vc)\in
\{0, +\infty\}.
\end{eqnarray}
In particular, one obtains $\eta(G)\in \{0, +\infty\}$.
\end{theorem}

Moreover, in Section \ref{s:random}, we prove some generic results for
random trees and their deterministic Schr\"odinger operators. We obtain:
\begin{proposition}\label{p:rangene}
  Let $G=(E,V)$ be a random tree with independent and identically
  distributed (i.i.d.) offspring.
  Suppose that the offspring distribution has finite expectation.
  Then for almost all trees, the Schr\"odinger operators
  $\Ac_G+\Vc$ and $\Delta_G+V$ are essentially self-adjoint on
  $\Cc_c(G)$, for all potentials $\Vc:V\to\R$. In particular, almost
  surely, one gets $\eta(G)=0$.
\end{proposition}
We refer to Section \ref{s:random} for definitions, a proof of this
result and also for Proposition \ref{p:birthrand}, which treats the case
of random offspring at a given generation.

We now present the structure of the paper. We start by proving, in
Section \ref{s:noneq}, that the domains of the Laplacian and of the
adjacency matrix are different for simple graphs of unbounded degree.
Then, in Section \ref{s:essSA}, we prove Proposition
\ref{p:essSA}. Next, we present our main tool in Section
\ref{s:surgery}. Subsequently, after giving a few definitions in
Section \ref{s:def}, we discuss the setting of trees. We start by
explaining in Section \ref{s:birth} how to reduce in some cases the
analysis of adjacency matrices to the one of Jacobi matrices. After
that, in Section \ref{s:expgro}, we provide an example of tree $G$
with ``exponential growth'' such that $\eta(G)=\infty$. Then, we prove
Proposition \ref{p:rangene} in Section \ref{s:random}. Next, in
Section \ref{s:possible}, we prove the first main result of the
introduction, namely Theorem \ref{t:main0} and generalize it in
Section \ref{s:recu}. Finally in Appendix~\ref{s:stabind}, we recall
a general result of stability of deficiency indices,
proposition~\ref{p:stab}, and deduce a criterion for essential
self-adjointness of Jacobi matrices, which possess unbounded
diagonals.

\noindent\textbf{Notation:}
The set of nonpositive integers is denoted by $\N$, note that  
$0\in\N$. Given a set~$X$ and $Y\subseteq X$ let $\bone_Y\colon X\to\{0,1\}$
be the characteristic function of~$Y$, namely $\bone_Y^{-1}\{1\}=Y$.
We denote also by $Y^c$ the complement set of~$Y$ in~$X$.

\vskip1mm
\noindent\textbf{Acknowledgments:} We would like to thank Hermann
Schulz-Baldes, Vladimir Georgescu, Andreas Knauf, and Mathias Rafler
for helpful discussions.

\section{General results}\label{s:adj}
\subsection{Comparison of domains}\label{s:noneq}
In view of Proposition \ref{p:stab}, it is tempting to try to prove that
the adjacency matrix~$\Ac_G$ is self-adjoint by comparing it to the
discrete Laplace operator $\Delta_G$. (Remember that the latter is
always essentially self-adjoint on $\Cc_c(G)$ by Proposition \ref{p:essSA}.)
But, as a matter of fact, if the graph~$G$ is simple
and has unbounded degree, we prove in this section that this is
impossible.

Given a locally finite graph $G=(E,V)$ and a potential $\Vc:V\to\R$,
we set $\Hc_G:=\Ac_G+\Vc$.
We first recall that the domain of the adjoint is given by
\begin{eqnarray*}
  \Dc\bigl((\Hc_G)^*\bigr)=\Big\{f\in \ell^2(G),
  x\mapsto \left(\sum_{y\sim x}E(x,y)f(y)\right)+\Vc(x)f(x)\in\ell^2(G)\Big\}.
\end{eqnarray*}
Then, given $f\in\Dc((\Hc_G)^*)$, one has:
\begin{eqnarray*}
  \left((\Hc_G)^*f\right)(x)=\left(\sum_{y\sim
      x}E(x,y)f(y)\right)+\Vc(x)f(x),\text{ for all }x\in V.
\end{eqnarray*}
We prove the result:
\begin{proposition}\label{p:noneq}
  Consider $G=(E,V)$ and suppose there is a sequence $(x_n)_{n\in \N}$
  of points in~$V$, so that
  \begin{eqnarray}\label{e:noneq0}
    \lim_{n\to\infty}\sum_{y\sim x_n}E^2(y,x_n)=\infty\text{ and }
    \lim_{n\to\infty}\frac{\left(\sum_{y\sim x_n}E(y,x_n)\right)^2}
                          {\sum_{y\sim x_n}E^2(y, x_n)}=\infty.
  \end{eqnarray}
  Then, $\Dc(\Delta_G)\neq\Dc(\Ac_G)$.
  In particular, the conclusion holds true when~$G$
  is simple and has unbounded degree.
\end{proposition}
\proof We suppose that $\Dc(\Delta_G)=\Dc(\Ac_G)$. Therefore, the uniform
boundedness principle ensures that there are $a,b\geq 0$ such that
\begin{eqnarray}\label{e:noneq1}
 \norm{\Delta_G f}^2\leq a\norm{\Ac_G f}^2 + b\norm{f}^2, \text{ for all }f\in
 \Dc(\Ac_G).
\end{eqnarray}
We note now that one has that
$\norm{\Delta_G(\bone_{\{x_n\}})}^2=\sum_{y\sim x_n}E^2(y,x_n)
  +\bigl(\sum_{y\sim x_n}E(y,x_n)\bigr)^2$
and also that $\norm{\Ac_G(\bone_{\{x_n\}})}^2= \sum_{y\sim x_n}
E^2(y,x_n)$. Taking $f=\bone_{\{x_n\}}$ in \eqref{e:noneq1} leads to a
contradiction.

Finally, when~$G$ is simple and has unbounded degree, consider a sequence
$(x_n)_{n\in \N}$, so that $\dG(x_n)$ tends to infinity. \qed

\subsection{Essential self-adjointness of discrete Schr\"odinger operators}
\label{s:essSA}
We prove some criteria of self-adjointness for Schr\"odinger operators.

\proof[Proof of Proposition \ref{p:essSA}]
We start with the first point and mimic \cite[Theorem
1.3.1]{Woj}. Using Proposition \ref{p:stab}, it is enough to suppose that
$\Vc$ is non negative. Take $f\in\Dc\big((\Delta_G+\Vc)^*\big)$ so
that $(\Delta_G+\Vc)^*f=-f$. Since $\Delta_G+\Vc$ is non-negative, it
is enough to prove that $f=0$. Notice that one has, for all $x\in V$,
\begin{eqnarray*}
  \sum_{y\sim x}E(x,y)f(y)=
  \Bigl(1+\Vc(x)+\sum_{y\sim x}E(x,y)\Bigr)f(x).
\end{eqnarray*}
Therefore, given $x\in V$, there exists $y\sim x$ with $|f(y)|>|f(x)|$.
This is in contradiction to the fact that $f\in\ell^2(G)$.

We turn to the second point. As there is no restriction on~$\Vc$,
it is enough to consider the case of $H:=\Ac_G+\Vc$.  We denote by
$\sphere i:=\{x\in\Vc,\dV(x_0,x)=i\}$ the sphere of radius $i\in\N$ around $x_0\in V$.
For $n\in\N$, consider $a_n:\N\to[0,1]$ with finite support and set
$\cchi_n:=\sum_{i\in \N} a_{n}(i)\bone_{\sphere i}$ and $\tilde\cchi_n:=1-\cchi_n$.
We see immediately that $\cchi_n\Dc(H^*)\subset\Dc(H)\subset\Dc(H^*)$.
Then, the commutator $[H^*,\cchi_n]$, defined on $\Dc(H^*)$,
is well defined (in the operator sense). Easily, it extends to a
bounded operator, which we denote by $[H^*,\cchi_n]_\circ$.
We take $f\in \Dc(H^*)$ and will prove that it is
also contained in $\Dc(H)$ by approximating it with $f_n:=\cchi_nf$.
We have
\begin{equation}\label{e:todo}
  \begin{split}
    \norm{f_m-f_n}+\norm{H(f_m-f_n)}\quad\leq&\quad\norm{(\cchi_m-\cchi_n) f}+\norm{(\cchi_m-\cchi_n)H^*f}\\&\quad{}
      +\norm{[H^*,\cchi_n]f}+\norm{[H^*,\cchi_m]f}\text.
  \end{split}
\end{equation}
We now choose $a_n$ in order to make $(f_n)_{n\in\N}$
a Cauchy sequence with respect to the graph norm of~$H$.
Set
\begin{equation}\label{e:an}
  a_n (i):=\begin{cases}
    1,& \mbox{ for } i \leq n,\\
    \min \bigl\{1, \max\{ 0, 1 - \frac1n\sum_{j=n+1}^i 1/b_j\}\bigr\},
      & \mbox{ for } i > n.
  \end{cases}
\end{equation}
Notice that $a_n$ has finite support, since $\sum_{j\in \N}1/b_j=+\infty$.
This gives that $\cchi_nf$ and $\cchi_nH^*f$
tend to~$f$ and $H^* f$, respectively.
It remains to control the commutator in \eqref{e:todo}.
By the Schur test and \eqref{e:an}, we have:
\begin{align*}
  \norm{[H^*, \tilde \cchi_n]_\circ}&\leq \sup_{v\in V} \sum_{w\in V}
  |\langle \bone_{\{v\}}, [H^*, \tilde \cchi_n]_\circ\,
  \bone_{\{w\}}\rangle|= \sup_{v\in V} \sum_{w\in V}
  |\langle \bone_{\{v\}}, [\Ac_G, \tilde \cchi_n]
  \bone_{\{w\}}\rangle|
  \\
  &= \sup_{v\in V} \sum_{w\in V} E(v, w) |\cchi_n(w)- \cchi_n(v)| =
  \sup_{v\in V} \sum_{w\in V, \dV(w, v)=1} E(v, w) |\cchi_n(w)-
  \cchi_n(v)|\leq\frac1n\text.
\end{align*}
Returning to \eqref{e:todo},
this implies that $f_n$ is a Cauchy sequence in $\Dc(H)$.
Let~$g$ be its limit.  Since~$H$ is closed, $g\in \Dc(H)$ and $g=f$.

We turn to (3) and (4). Taking in account the contribution of
the potential, we essentially rewrite \cite[Proposition 1.1]{Gol}.
Take $f\in \Cc_c(G)$. For $\dG$ bounded let $\Mr(x):=E(x)$
and for~$E$ bounded let $\Mr(x):=\dG(x)$.
Let $\Mr$ be the operator of multiplication by $\Mr(\cdot)$, too.
We denote all constants, which are independent from~$f$, by the same letter~$C$.
We have:
\begin{align*}
  \norm{(\Ac_G+\Vc)f}^2&
  \le2\sum_x\Bigl|\sum_{y\sim x}E(x,y)f(y)\Bigr|^2+2\norm{\Vc f}^2
  \le2\sum_x\dG(x)E^2(x)\sum_{y\sim x}|f(y)|^2+2\norm{\Vc f}^2\\&
  \le2\sum_x\dG(x)\max_{y\sim x}(\dG(y))E^2(x)|f(x)|^2+2 \norm{\Vc f}^2\\&
  \le2\sum_xE^2(x)\dG(x)\bigl(C+\dG(x)\bigr)|f(x)|^2+2\norm{\Vc f}^2
  \le C\norm{\Mr f}^2.
\end{align*}
Moreover, noticing that the potential~$\Vc$ commutes with~$\Mr$, we get
\begin{align*}
  |\langle f,[\Ac_G,\Mr]f\rangle|&
    =\Bigl|\sum_x\overline{f(x)}
           \sum_{y\sim x}E(x,y)\big(\Mr(y)-\Mr(x)\big)f(y)\Bigr|
  \le\sum_x\sum_{y\sim x}C|E^{1/2}(x)f(x)|\,|E^{1/2}(y)f(y)|\\&
  \le c\sum_x\dG(x)|E^{1/2}(x)f(x)|^2
  \leq C\bignorm{\Mr^{1/2}f}^2.
\end{align*}
Then, using \cite[Theorem X.36]{RS}, the result follows.

We deal now with the fifth point. As a potential is essentially
self-adjoint on $\Cc_c(G)$, thanks to W\"ust's Theorem, e.g.,
\cite[Theorem X.14]{RS}, it is enough to prove that there is $b\ge0$
so that,
\begin{eqnarray}\label{e:Wuest}
  \norm{\Ac_Gf}^2\le\norm{\Vc f}^2+b\norm{f}^2, \mbox{ for all } f\in \Cc_c(G).
\end{eqnarray}
As $x\mapsto\Vc(x)\bone_{K}(x)$ is bounded, it is enough to prove
\eqref{e:Wuest} with $b=0$ and under the stronger hypothesis:
$\sum_{y\sim x}E^2(x,y)\dG(y)\le\Vc^2(x)$ for all $x\in V$.
The statement is now obvious as, for all $f\in\Cc_c(G)$, one has
\begin{eqnarray*}
  \norm{\Ac_G f}^2
    =\sum_{x\in V}\sum_{y\sim x}|E(x,y)f(y)|^2
    \le\sum_{x\in V}\sum_{y\sim x}\dG(x)E^2(x,y)|f(y)|^2
    =\sum_{x\in V}\sum_{y\sim x}\dG(y)E^2(x,y)|f(x)|^2.
\end{eqnarray*}
Finally, by using the last inequality and by taking into account
the diagonal part of the Laplacian, one has, for all $f\in\Cc_c(G)$,
\begin{eqnarray*}
  \norm{\Delta_G f}^2
    \leq\sum_{x\in V}\Bigl(\sum_{y\sim x}\dG(y)E^2(x,y)+E^2(x,y)\Bigr)|f(x)|^2.
\end{eqnarray*}
W\"ust's Theorem gives the last point.
\qed

\begin{remark}
Given $a\in [0,1)$, note that if one strengthens the assumption
in the fourth point to
\[\sum_{y\sim x}E^2(x,y)\dG(y)\le a\Vc^2(x),\mbox{ for all }x\notin K\]
the previous proof and the Kato-Rellich theorem
(or more generally Proposition \ref{p:stab}) ensures
$\Dc(\Ac_G+\Vc)=\Dc(V)$, too.
In the same spirit, if one supposes that
$\sum_{y\sim x}E^2(x,y)\big(1+\dG(y)\big)\leq a\Vc^2(x)$
for all $x\notin K$ in the fifth point,
one gets also $\Dc(\Delta_G+\Vc)=\Dc(V)$.
\end{remark}

\subsection{Bounded perturbations of graphs and deficiency
  indices}\label{s:surgery}
In this section, we compute the deficiency indices,
in the case one adds up to a given number of edges per vertex
to a countable union of graphs.
We slightly improve the surgery Lemma of \cite{Gol}.

\begin{lemma}\label{l:surgery}
  Given a sequence of graphs $G_n=(E_n, V_n)$, for $n\in \N$, let
  $G^\circ:=(E^\circ,V^\circ):=\bigcup_{n\in\N}\,G_n$
  be the disjoint union of $\{G_n\mid n\in\N\}$.
  Choose $\tilde E:V^\circ\times V^\circ\to[0,\infty)$,
  so that $\tilde E$ is symmetric, with support away from the diagonal.
  Set $G:=(E,V)$ with $V=V^\circ$ and $E:=E^\circ+\tilde E$.
  Suppose that:
  \begin{eqnarray}\label{e:surgery}
    \sup_{x\in V}\sum_{y\in V}\tdG(y){\tilde E}^2(x,y)<\infty,
  \end{eqnarray}
  where $\tdG(x):=|\{y\in V,\tilde E(x,y)\ne0\}|$.
  Consider a potential $\Vc:V\to\R$.
  Set $\Hc_G:=\Ac_G+\Vc$ and $\Hc_{G_n}:=\Ac_{G_n}+\Vc|_{G_n}$.
  Then, one obtains
  \begin{eqnarray*}
    \eta(\Hc_G)=\sum_{n\in\N}\eta(\Hc_{G_n}).
  \end{eqnarray*}
  In particular, $\eta(G)=\sum_{n\in\N}\eta({G_n})$.
\end{lemma}
\proof Take $f\in\Cc_c(G)=\Cc_c(G^\circ)$.
Set $\Hc_{G^\circ}:=\bigoplus_{n\in\N}\Hc_{G_n}$.
Notice that:
\begin{align*}
  \norm{(\Hc_{G}-\Hc_{G^\circ})f }^2&
    =\sum_{x\in V}\Bigl|\sum_{y\in V}\tilde E(x,y)f(y)\Bigr|^2
    \le\sum_{x\in V}\sum_{y\in V}\tdG(x)\tilde E^2(x,y)|f(y)|^2\\&
    =\sum_{x\in V}\Bigl(\sum_{y\in V}\tdG(y)\tilde E^2(x,y)\Bigr)|f(x)|^2.
\end{align*}
We infer, there is a finite~$M$,
so that $\norm{(\Hc_G-\Hc_{G^\circ})f}\leq M\norm{f}$,
for all $f\in\Cc_c(G)=\Cc_c(G^\circ)$.
Then, the closure of $(\Hc_G-\Hc_{G^\circ})$
is a bounded operator and Proposition~\ref{p:stab} can be applied.

Alternatively, one can conclude using an argument of \cite{Gol}.
Since the closure of $(\Hc_G-\Hc_{G^\circ})$
is a bounded operator, the graph norms of $\Hc_G$ and of
$\Hc_{G^\circ}$ are equivalent when restricted to $\Cc_c(G)$.
By taking the closure, we infer $\Dc(\Hc_G)=
\Dc(\Hc_{G^\circ})$. Moreover, using again the boundedness of the
difference and the definition of the domain of the adjoints of $\Hc_G$
and of $\Hc_{G^\circ}$, one gets directly $\Dc((\Hc_G)^*)=
\Dc((\Hc_{G^\circ})^*)$. Finally, since the deficiency indices
$\eta_{\pm}(\Hc_G)$ of $\Hc_G$ are equal (and of $\Hc_{G^\circ}$,
resp.), \eqref{e:dom} gives that $\eta(\Hc_G)= \eta(\Hc_{G^\circ})$.
\qed

\begin{example}\label{ex:remove}
  Given a locally finite graph $G:=(E,V)$ with bounded weights~$E$
  and a set of vertices $X\subseteq V$, such that $\sup\dG(X)<\infty$,
  then the induced graph $G'=G[V\setminus X]$,
  obtained by removing the vertices in~$X$,
  has deficiency index $\eta(G')=\eta(G)$.
\end{example}

\subsection{Tensor products and deficiency indices}
For the sake of completeness and motivated by Remark~\ref{r:tensor}
(see below), we discuss shortly the tensor product of graphs regarding
the computation of deficiency indices. We
recall that given two graphs
$G_i=(E_i, V_i)$, $i=1,2$, one defines the tensor product $G:=(E,V)$
of $G_1$ with $G_2$ by setting $V:=V_1\times V_2$ and $E\big((x_1,
x_2), (y_1, y_2)\big):= E(x_1, y_1)\cdot E(x_2, y_2)$. One sees that
$\Ac_{G_1\otimes G_2}= \Ac_{G_1}\otimes \Ac_{G_2}$. We turn to the
question of deficiency indices. It is well-known that $\eta(G_1\otimes
G_2)=0$ if $\eta(G_1)=\eta(G_2)=0$, e.g., \cite[Theorem
VII.33]{RS}. One has also that $\eta(G_1\otimes G_2)= \infty$, when
$\eta(G_1)= \infty$ and $\eta(G_2)>0$. In fact, in the general case,
one obtains easily a lower bound on the deficiency indices:

\begin{lemma}\label{l:known}
  Given two symmetric operators $S,T$ acting on the Hilbert spaces~$\Hr$
  and $\Kr$, respectively.
  Let $\eta=\max_{i\in\{\pm\}}\bigl(\eta_i(S)\cdot\eta_i(T)\bigr)$,
  with the convention $0\cdot\infty=0$.
  Then, $\eta_{\pm}(S\otimes T)\geq\eta$.
\end{lemma}
\proof We recall that, given a symmetric operator~$H$,
$z\mapsto\dim\ker(H^*-z)$ is
constant on the upper and lower open half-planes of $\C$.
Therefore it is enough to give a lower bound for
$\dim\ker(S^*\otimes T^*-z^2)$, for $z= e^{\rmi\pi(1/2\pm1/4))}$.
Take $f\in \Dc(S^*)$ and $g\in \Dc(T^*)$, so
that $S^*f=zf$ and $T^*g=zg$. One has:
\begin{eqnarray*}
S^*\otimes T^* (f\otimes g)- z^2 f\otimes g = (S^*f-zf)\otimes T^*g + z
f\otimes (T^*g-zg)=0.
\end{eqnarray*}
This concludes the proof. \qed

It is however more important to obtain the exact value of the
deficiency indices. We recall the following elementary fact:

\begin{lemma}\label{l:tensor}
  Let~$G$ be a locally finite graph and~$K$ be a finite graph.
  Then, one deduces
  \begin{equation*}
    \eta({G\otimes K})=\eta(G)\cdot\dim(\mathrm{Im}(\Ac_K)).
  \end{equation*}
\end{lemma}
\proof As $\Ac_K$ is self-adjoint in a finite dimensional Hilbert space,
we can decompose it with the help of its eigenspaces.
We have $\Ac_K=\bigoplus_i\lambda_i\bone_{E_i}$,
where $E_i$ is the eigenspace associated to the eigenvalue~$\lambda_i$.
Note that $(\Ac_{G\otimes K})^*=\bigoplus_i\lambda_i(\Ac_G)^*\otimes\bone_{E_i}$.
To conclude, we notice that
$\dim\big(\ker((\Ac_G)^*\otimes\bone_{E_i}+\rmi)\big)
  =\dim\big(\ker((\Ac_G)^*+\rmi)\big)\times \dim \bone_{E_i}$.
\qed

We now come back to the conjecture mentioned in the introduction
following \eqref{e:problem}

\begin{remark}\label{r:tensor}
  The complete graph $K_2=(E_2,V_2)$ is defined by $V_2:=\{0,1\}$
  and $E_2(0,1)=1$.  Note that $\Ac_{K_2}$ is injective.
  Its spectrum is $\{-1,1\}$.
  Given a locally finite graph~$G$, the previous lemma
  states that $\eta(G\otimes K_2)=2\eta(G)$.
  Moreover, note that $G\otimes K_2$ is bipartite.
  Therefore if \eqref{e:problem} is true for all bipartite simple graphs,
  then it is true for all simple graphs.
\end{remark}

\section{The case of a tree}
\subsection{Some definitions related to trees}\label{s:def}
It is convenient to choose a root in the tree.
Due to its structure, one can take any point of~$V$. We denote it by $\epsilon$.

We define inductively the \emph{spheres} $\sphere n$ by $\sphere{-1}=\emptyset$,
$\sphere0:=\{\epsilon\}$, and $\sphere{n+1}:=\Nr_G(\sphere n)\setminus\sphere{n-1}$.
Given $n\in\N$, $x\in\sphere n$, and $y\in\Nr_G(x)$,
one sees that $y\in\sphere{n-1}\cup\sphere{n+1}$.
We write $x\sim>y$ and say that~$x$ is a \emph{son} of~$y$, if
$y\in\sphere{n-1}$, 
while we write $x<\sim y$ and say that~$x$ is a \emph{father} of~$y$,
if $y\in\sphere{n+1}$. 
Notice that $\epsilon$ has no father.
Given $x\ne\epsilon$, note that there is a unique $y\in V$ with $x\sim>y$,
i.e., everyone apart from $\epsilon$ has one and only one father.
We denote the father of~$x$ by $\father x$.
Given $x\in\sphere n$, we set $\ell(x):=n$, the \emph{length} of~$x$.
The \emph{offspring} of an element~$x$ is given by
$\off(x):=|\{y\in \Nr_G(x),y\sim>x\}|$, i.e., it is the number of sons
of~$x$. When $\ell(x)\ge1$, note that $\off(x)=\dG(x)-1$. 

\subsection{Diagonalization in the case of an offspring
  depending on the generation}\label{s:birth}
In this section, we define a certain family of trees.
Then, we explain how to explicitly diagonalize the adjacency matrices on them.
We start with a definition.

\begin{definition}\label{d:birth}
  A simple tree $G=(E,V)$ with \emph{offspring sequence}~$(b_n)_{n\in\N}$
  is a simple tree with a root such that $b_n=\off(x)$,
  for each $x\in\sphere n$ and $n\in\N$.
\end{definition}
In Proposition \ref{p:birthrand}, we consider a family of trees
with random offspring per individual and generation.
At the moment, we focus on the deterministic case and give a concrete example:
\begin{align*}
  \xymatrix{%
    &&\epsilon\ar@{-}[d]\ar@{-}[dr]&&&& S_0
    \\
    &&{\bullet}\ar@{-}[d]\ar@{-}[dl]\ar@{-}[dll]&
      {\bullet}\ar@{-}[rrd]\ar@{-}[d]\ar@{-}[dr]&&& S_1
    \\
   {\bullet}\ar@{.}[d]&
   {\bullet} \ar@{.}[d]&
   {\bullet}\ar@{.}[d]&
   {\bullet}\ar@{.}[d]&
  {\bullet} \ar@{.}[d]&
  {\bullet} \ar@{.}[d]& S_2
    \\
    &&&&&&
  }
  \\
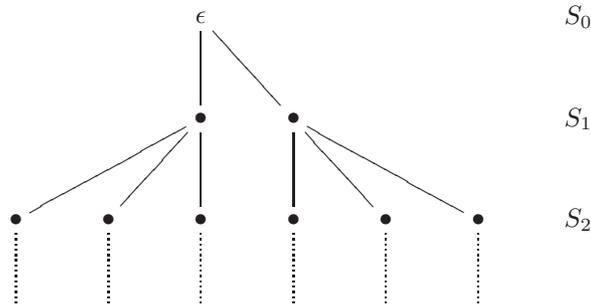
  \text{\it Example of a tree with } b_0=2 \text{\it\, and } b_1=3.
  \quad\quad
\end{align*}
Now we adapt the decomposition of a tree given in \cite{AF}, see also \cite{GG},
in order to write the adjacency matrix as a direct sum of Jacobi matrices.
We consider the tree $G=(E,V)$ with offspring sequence~$(b_n)_{n\in\N}$.
We define:
\begin{eqnarray*}
  (U f)(x)
    :=\bone_{\{\cup_{n\ge1}\sphere n\}}(x)\frac1{\sqrt{b_{\ell\big(\father
        x\big)}}}\,f\big(\father x\big),
    \text{ for }f\in\ell^2(G).
\end{eqnarray*}
Easily, one get $\norm{U f}=\norm f$, for all
$f\in\ell^2(G)$. Moreover, it is a completely non-unitary isometry, i.e., it
is an isometry, such that the strong limit
$\slim\limits_{k\to\infty}(U^*)^k=0$. The adjoint~$U^*$ of~$U$ is
given by  
\begin{eqnarray*}
  (U^*f)(x):=\frac1{\sqrt{b_{\ell(x)}}}\sum_{y\sim>x}f(y),
    \text{ for }f\in\ell^2(G).
\end{eqnarray*}
Note that one has:
\begin{eqnarray*}
  (\Ac_Gf)(x)=\sqrt{b_{\ell(\father x)}}\,(Uf)(x)+\sqrt{b_{\ell(x)}}\,(U^*f)(x),
    \text{ for }f\in\Cc_c(G).
\end{eqnarray*}
Supposing now that $b_n\geq 1$ for all $n\in \N$,
we construct invariant subspaces for $\Ac_G$.
We start by noticing that $\dim \ell^2(S_n)=\prod_{i=0, \ldots, n-1} b_n$,
for $n\geq 1$ and $\dim \ell^2(S_0)=1$.  Therefore, as $U$ is an isometry,
$U\ell^2(S_n)=\ell^2(S_{n+1})$ if and only if $b_n=1$.
Set $\Q_{0,0}:=\ell^2(\sphere0)$ and $\Q_{0,k}:=U^k\Q_{0,0}$, for all $k\in\N$.
Note that $\dim \Q_{0,k}= \dim \ell^2(\sphere0)=1$, for all $k\in\N$.
Moreover, given $f\in \ell^2(S_k)$, one has $f\in \Q_{0,k}$ if and
only if $f$ is constant on $S_k$. We define recursively $\Q_{n,n+k}$
for $k,n\in\N$.  Given $n\in\N$, suppose that $\Q_{n,n+k}$ is
constructed for all $k\in \N$, and set
\begin{itemize}
  \item $\Q_{n+1,n+1}$ as the orthogonal complement of
    $\bigoplus_{i=0,\ldots,n}\Q_{i,n+1}$ in $\ell^2(\sphere{n+1})$,
  \item $\Q_{n+1,n+k+1}:= U^k\Q_{n+1,n+1}$, for all $k\in \N\setminus\{0\}$.
\end{itemize}
We sum-up the construction in the following diagram:
\vspace*{-.7cm}
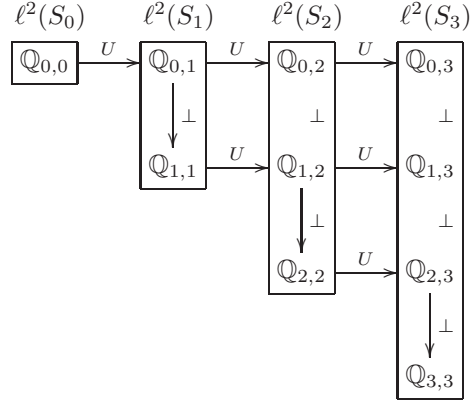
\begin{align*}
\def\gvide#1{\save
[].[]!C="gvide#1"*[F]\frm{}\restore}
\def\g#1{\save
[].[d]!C="g#1"*[F]\frm{}\restore}%
\def\gg#1{\save
[].[dd]!C="gg#1"*[F]\frm{}\restore}
\def\ggg#1{\save
[].[ddd]!C="ggg#1"*[F]\frm{}\restore}
\xymatrix{%
\ar@{}[d]^{\hspace*{-0.5cm}\txt{$\ell^2(S_0)$}}&
\ar@{}[d]^{\hspace*{-0.5cm}\txt{$\ell^2(S_1)$}}
&\ar@{}[d]^{\hspace*{-0.5cm}\txt{$\ell^2(S_2)$}}&
\ar@{}[d]^{\hspace*{-0.5cm}\txt{$\ell^2(S_3)$}}
\\
\gvide1 \Q_{0,0}\ar@{->}[r]^{U}& \g1\Q_{0,1}\ar@{->}[r]^{U}\ar@{->}[d]^{\perp}&
\gg2 \Q_{0,2}\ar@{->}[r]^{U}\ar@{}[d]^{\perp}&\ggg1\Q_{0,3} \ar@{}[d]^{\perp}
\\
%\ell^2(S_0)\ar@{-->}[r]
& \Q_{1,1}\ar@{->}[r]^{U}&\Q_{1,2}\ar@{->}[r]^{U}\ar@{->}[d]^{\perp}&\Q_{1,3}\ar@{}[d]^{\perp}
\\
&& \Q_{2,2}\ar@{->}[r]^{U}&\Q_{2,3} \ar@{->}[d]^{\perp}
\\
&& &\Q_{3,3}}
\end{align*}

We point out that $\dim \Q_{n+1,n+1}= \dim \Q_{n+1,n+k+1}$, for all
$k\in \N$ and stress that it is $0$ if and only if $b_n=1$.
Notice that $U^*\Q_{n,n}=0$, for all $n\in\N$.
Set finally $\M_n:=\bigoplus_{k\in\N}\Q_{n,n+k}$
and note that $\ell^2(G)=\bigoplus_{n\in\N}\M_n$.
Moreover, one has that canonically
$\M_n\simeq\ell^2(\N;\Q_{n,n})\simeq\ell^2(\N)\otimes\Q_{n,n}$.
In this representation, the restriction~$\Ac_n$ of $\Ac_G$
to the space~$\M_n$ is given by the following tensor product of Jacobi matrices:
\begin{eqnarray*}
  \Ac_n\simeq
    \left(\begin{array}{cccccc}
    0 & \displaystyle \sqrt{b_{n}}& 0 & 0 &\cdots
    \\
    \displaystyle \sqrt{b_{n}}& 0 &\displaystyle
    \sqrt{b_{n+1}} & 0& \ddots
    \\
    0 & \displaystyle\sqrt{b_{n+1}}& 0 & \displaystyle
    \sqrt{b_{n+2}}& \ddots
    \\
    \vdots & \ddots &\ddots & \ddots & \ddots
    \end{array}\right) \otimes \bone_{\Q_{n,n}}.
\end{eqnarray*}
Now $\Ac_G$ is given as the direct sum $\bigoplus_{n\in\N}\Ac_n$
in $\oplus_{n\in\N}\M_n$.
In particular, $\eta(G)=\sum_{n\in\N}\eta(\Ac_n)$. Note that if
  $b_n=1$, for all $n\in \N$, we recover the case of the adjacency
  matrix of the simple graph $\N$.

We now turn to the case of $b_n:=\lfloor n^\alpha\rfloor$,
for some $\alpha>0$.

\proof[Proof of Proposition~\ref{p:tree}]
The sum $\sum_{n\in\N}\ 1/\sqrt{b_n}$ is finite if and only if
$\alpha>2$. Then \cite[page 504]{Be} yields that $\Ac_n=(\Ac_n)^*$ for
$\alpha\in[0,2]$ and $n\in\N$. One infers that $\eta(G)=0$.
Now, easily one sees that $b_{i-1}b_{i+1}\le b_i^2$, for $i\ge1$.
Thus, \cite[page 507]{Be} gives that $\eta(\Ac_n)= \dim(\Q_{n,n})$.
This completes the proof. \qed

\subsection{Trees with exponential growth and non-essential
  self-adjointness}\label{s:expgro}

In the previous section, we focused on trees with given offspring per
individual for each generation. We now replace this hypothesis by a
control on the maximum and on the minimum of the offspring of
individuals for each generation. We turn to the result, see also \cite{MO, Mu}.

\begin{proposition}\label{p:exp}
  Let $G=(E,V)$ be a locally finite simple tree, endowed with an origin.
  Supposing
  \begin{eqnarray}\label{e:exp}
    n\mapsto\frac{\max_{x\in\sphere{n-1}}\off(x)}{\min_{x\in\sphere n}\off(x)}\in\ell^1(\N),
  \end{eqnarray}
  one has that $\eta(G)=\infty$.
\end{proposition}
Condition \eqref{e:exp} can be interpreted as an ``exponential growth''.

\proof We construct $f\in\ell^2(G)\setminus \{0\}$,
so that $(\Ac_G)^*f=\rmi f$ and
\begin{eqnarray*}
  f(x)=f(y),\text{ if }\father x=\father y,\quad\quad(x,y\in
  V\setminus\{\epsilon\})
\end{eqnarray*}
i.e., for all $x\in V$, $f$ is constant on $\off(x)$.
We denote the constant value by $f(\sim>x)$.
With this nota\-tion, we have
\begin{eqnarray}\label{e:const0}
  \off(x)f(\sim>x)+f\big(\father x\big)=\rmi f(x),
\end{eqnarray}
for all $x\in\sphere n$, with $n\ge1$.
We denote by $\norm[\sphere n]{f}$ the $\ell^2$-norm of~$f$ restricted to~$\sphere n$.
Then we have:
\begin{align*}
  \norm{f}^2_{\sphere{n+1}}&
    =\sum_{x\in\sphere{n-1}}\sum_{y\sim>x}\sum_{z\sim>y}|f(z)|^2
    =\sum_{x\in\sphere{n-1}}\sum_{y\sim>x}\sum_{z\sim>y}|f(\sim>y)|^2,\\&
    \le\sum_{x\in\sphere{n-1}}\sum_{y\sim>x}\sum_{z\sim>y}
      \frac2{\off^2(y)}\left(|f(y)|^2+|f(x)|^2\right),
        \text{ by \eqref{e:const0}}\\&
    \le2\frac{\max_{x\in\sphere{n-1}}\off(x)}{\min_{x\in\sphere n}\off(x)}
      \norm{f}^2_{\sphere{n-1}}+\frac2{\min_{x\in\sphere n}\off(x)}\norm{f}^2_{\sphere n}.
\end{align*}
By induction, one sees that $\sup_{n\in \N}\norm[\sphere n]{f}^2$ is finite.
Finally using \eqref{e:exp}, we derive that $f\in\ell^2(G)$.
Theorem~\ref{t:main0} concludes that the deficiency indices are infinite.
\qed

\subsection{Discrete Schr\"odinger operators and random trees}\label{s:random}

In this section we discuss certain random trees.
Before dealing with random trees in the sense of Definition~\ref{d:randomtree},
we start with trees with random offspring sequence,
see Definition~\ref{d:birth}.

We recall some well-known notions from probability theory.
The \emph{left shift} on $\N^\N$ is $\tau\colon\N^\N\to\N^\N$,
$\tau\bigl((x_n)_{n\in\N}\bigr):=(x_{n+1})_{n\in\N}$.
We assign the discrete topology to~$\N$ and the product topology
to~$\N^\N$. Therefore, $\tau$ is continuous. An~$\N$-valued stochastic
process $X:=(X_n)_{n\in\N}$, is called \emph{ergodic}, if for all
Borel-measurable $A\subseteq\N^\N$, one has
\begin{align*}
  P\bigl(X\in A\text{ and }\tau(X)\not\in A\bigr)
      +P\bigl(X\not\in A\text{ and }\tau(X)\in A\bigr)
    =0 \implies P(X\in A)\in\{0,1\}
\end{align*}
and \emph{stationary}, if
\begin{equation*}
  P(X\in A)=P\bigl(\tau(X)\in A\bigr)
\end{equation*}
for all Borel-measurable $A\subseteq\N^\N$.
For example, if $X_n$, $n\in\N$, are i.i.d.\ random variables
then the process $(X_n)_{n\in\N}$ is stationary and ergodic.

\begin{proposition}\label{p:birthrand}
  Let $G=(E,V)$ be a tree with offspring sequence~$(b_n)_{n\in\N}$,
  where $(b_n)_{n\in\N}$ is a stationary and ergodic stochastic process.
  Then for almost every~$G$, the Schr\"odinger operators
  $\Ac_G+\Vc$ and $\Delta_G+\Vc$ are essentially self-adjoint on
  $\Cc_c(G)$, for all $\Vc:V\to \R$.
\end{proposition}
\proof
%Given a stationary and ergodic $\N$-valued
%stochastic process $X:=(X_n)_{n\in\N}$ and $x\in\N$ with $P(X_0=x)>0$,
%there exists almost surely a subsequence of $(X_n)_{n\in\N}$
%which is constant with value~$x$.
%Indeed, let
%$M_x:=\{(x_n)_{n\in\N}\mid\forall n\in\N\colon x_n\ne x\}$ be the set
%of all sequences which never take the value~$x$. Then
%$M_x\subseteq\tau^{-1}(M_x)$ and
%$M_x\setminus\tau^{-1}(M_x)=\emptyset$. By stationarity, we see
%\begin{equation*}
%P\bigl(X\in\tau^{-1}(M_x)\setminus M_x\bigr)
%  =P\bigl(X\in\tau^{-1}(M_x)\bigr)-P\bigl(X\in M_x\bigr)
%  =0.
%\end{equation*}
%Hence ergodicity yields $P(X\in M_x)\in\{0,1\}$.
%Remembering $P(X_0=x)>0$ we entail $P(X\in M_x)=0$.
%By stationarity we infer
%\begin{equation*}
 % P\Bigl(X\in\bigcup\nolimits_{k\in\N}\tau^{-k}(M_x)\Bigr)
 % \le\sum\nolimits_{k\in\N}P\bigl(X\in\tau^{-k}(M_x)\bigr)
 % =\sum\nolimits_{k\in\N}P\bigl(\tau^k(X)\in M_x\bigr)
 % =0.
%\end{equation*}
%This implies that, for almost every sequence, the value~$x$
%is taken infinitely many times.
Take $m\in\N$, so that $P(b_0=m)>0$. Since $(b_n)_{n\in\N}$ is
a stationary and ergodic $\N$-valued stochastic process,
there is, almost surely, a subsequence $(b_{n_k})_{k\in\N}$
with $b_{n_k}=m$ for all $k\in\N$.
Consider now the forest of finite trees obtained by removing all edges
between $\sphere{n_k}$ and $\sphere{n_k+1}$, for all $n\in\N$.
Note that, for each element of~$\sphere{n_k+1}$,
there is at most one edge connecting it to $\sphere{n_k}$.
The Schr\"odinger operators, restricted to the finite trees,
are all essentially self-adjoint.
Lemma~\ref{l:surgery} gives the result.\qed

Next we consider random trees.
Denote by $W:=\bigcup_{n\in\N}(\N^*)^n$
the set of all finite words over the alphabet~$\N^*:=\N\setminus\{0\}$.
The~length of a word~$w=(w_1,\dotsc,w_n)\in W$ is $\ell(w):=n$.
\begin{definition}\label{d:randomtree}
  Let~$(X_w)_{w\in W}$ be a family of i.i.d.\ random variables
  with values in~$\N$.
  We construct a graph $G=(E,V)$ as follows:
  \begin{gather*}
    V:=\bigl\{(w_1,\dotsc,w_n)\in W\mid w_{m+1}\le X_{(w_1,\dotsc,w_m)}
      \text{ for all }m\in\N,m<n
    \bigr\}\text{ and}\\
    E(v,w):=\begin{cases}
      1&\text{if $\{\ell(v),\ell(w)\}=\{n,n+1\}$ and
	      $(v_0,\dotsc,v_n)=(w_0,\dotsc,w_n)$}\\
      0&\text{otherwise,}
    \end{cases}
  \end{gather*}
  for $v=(v_1,\dotsc,v_{\ell(v)}),w=(w_1,\dotsc,w_{\ell(w)})\in V$.
  We call~$G$ \emph{random tree with i.i.d.\ offspring}.
  The law of $X_\epsilon$ is called \emph{offspring distribution} of~$G$.
\end{definition}
Note that a random tree is a tree with the empty word~$\epsilon$ as root.
Words of length~$n$ correspond to $\sphere n$, the~$n$-sphere.
Hence, the notation $\ell$ of the length is consistent
with the one given in Section~\ref{s:def}.

\begin{proposition}\label{p:expect}
  Let $G=(E,V)$ be a random tree with i.i.d.\ offspring,
  such that its offspring distribution has finite expectation.
  Then almost surely there are $M\geq 1$ and a family $(G_i)_{i\in\N}$
  of disjoint finite subtrees $G_i:=(E_i,V_i)$ of~$G$, so that
  $V=\bigcup_{i\in\N}V_i$
\begin{eqnarray}\label{e:expect}
\sup_{i\in V_i}\max_{x\in \max(V_i)} {\rm off} (x) \leq M,
\end{eqnarray}
where $\max(V_i):=\{x\in V_i, (y \sim> x $ in $G) \implies y\notin V_i\}$,
for all $i\in \N$.
\end{proposition}
\proof
Since the offspring distribution has finite expectation,
there is $M\in\N$ such that
\begin{equation}\label{e:children}
  \sum\nolimits_{m>M}mP(X_\epsilon=m)<1.
\end{equation}
Let $\tilde G:=G\setminus L$
be the forest one gets by deleting all the edges in
\begin{equation*}
  L:=\{(v,w)\in V\times V\mid
    \ell(v)<\ell(w),\off_G(v)\le M\}
\end{equation*}
from~$G$.
Each connected component in $\tilde G$
is a random tree with independent offspring.
Denote by $\hat G=(\hat E,\hat V)$ a connected component of~$\tilde G$.
The expected number of sons in~$\hat G$
is given by the l.h.s.\ in~\eqref{e:children}.
It is well known that such family trees almost
surely get extinct, see e.g., \cite[Theorem 3.11]{Kle}.
Therefore all the connected components of $\tilde G$ are almost surely finite.
We present a proof here.

The tree~$\hat G$ has a root $\hat w_0\in\hat V$
with $\ell(\hat w_0)=\min\{\ell(\hat w)\mid\hat w\in\hat V\}$.
We define the~$n$-sphere of $\hat G$ to be
$\hatsphere n:=\{\hat w\in\hat V\mid\ell(\hat w)=n+\ell(\hat w_0)\}$
and denote by $\hat X_{\hat w}:=\off_{\hat G}(\hat w)$
the number of sons of $\hat w\in\hat V$ in~$\hat G$.
The random variable $\hat Y_n:=\lvert\hatsphere n\rvert$
fulfills $\hat Y_n=\sum_{\hat w\in\hatsphere{n-1}}\hat X_{\hat w}$
and is hence measurable with respect to the $\sigma$-algebra
$\hat\Fc_n:=\sigma\bigl(\hat X_{\hat w}\bigm|
\hat w\in\bigcup_{j=0}^{n-1}\hatsphere j\bigl)$.
Therefore the stochastic process $(\hat Y_n)_{n\in\N}$
is adapted to the filtration~$(\hat\Fc_n)_{n\in\N}$.
With \eqref{e:children}, for all $n\in\N$ we have
\begin{equation}\label{e:condexp}
\E\bigl[\hat Y_{n+1}\bigm|\hat\Fc_n\bigr]
  =\sum\nolimits_{\hat w\in\hatsphere n}\E\bigl[\hat X_{\hat w}
    \bigm|\hat\Fc_n\bigr]
  =\sum\nolimits_{\hat w\in\hatsphere n}\E[\hat X_{\hat w}]
  =\hat Y_n\E[\hat X_{\hat w_0}]\le\hat Y_n.
\end{equation}
Hence, the process $(\hat Y_n)_{n\in\N}$ is a supermartingale.
Since $\hat Y_n\ge0$, the martingale convergence
theorem guarantees that $Y_n$ converges almost surely.
We denote its limit by $\hat Y$.
With \eqref{e:condexp} we entail
\begin{equation*}
  0\leq\E[\hat Y_n]
  =\E\bigl[\E[\hat Y_n\mid\hat\Fc_{n-1}]\bigr]
  =\E\bigl[\hat Y_{n-1}\E[\hat X_{\hat w_0}]\bigr]
  =\E[\hat Y_{n-1}]\E[\hat X_{\hat w_0}]
  =\bigl(\E[\hat X_{\hat w_0}]\bigr)^n.
\end{equation*}
In view of \eqref{e:children}, Fatou's Lemma ensures that
$\E[\hat{Y}]=0$ and therefore that $\hat Y=0$ almost surely.
Finally, since $\hat Y$ assumes only integer values,
for almost every realization of~$(\hat Y_n)_{n\in\N}$
there exists $N\in\N$ with $\hat Y_n=0$ for all $n\ge N$.
\qed

It remains to prove the announced result.
\proof[Proof of Proposition \ref{p:rangene}]
Almost surely, Proposition \ref{p:expect} gives a
forest of finite trees $G_i=(E_i, V_i)$. On each of them, the
restriction of the Schr\"odinger operator is essentially self-adjoint,
as $\ell^2(G_i)$ is finite dimensional.
Moreover, as $\bigcup_{i\in\N}V_i=V$ and \eqref{e:expect} holds true,
the hypothesis of Lemma~\ref{l:surgery} are satisfied and the result
follows.\qed

\subsection{The possible indices}\label{s:possible}

We now prove our main result in the context of trees and improve it in
Section \ref{s:recu}. This is a proof by contradiction.

We start with some notations about subgraphs. The connected
component $\cc_G(x)$ of $x\in V$ is the graph $\cc_G(x):=(E_x,V_x)$
with $V_x:=\{y\in V,\text{ there is an $x$-$y$-path}\}$ and
$E_x:=E|_{V_x\times V_x}$. A graph $G':=(E',V')$ is called a
\emph{subgraph} of~$G$, if $V'\subseteq V$ and
$E'(x,y)\in\{0,E(x,y)\}$, for all $x,y\in V'$. The
subgraph $G[V']:=(E|_{V'\times V'},V')$ is called the \emph{induced
graph} of~$G$ by $V'\subseteq V$. Given a set $S\subseteq V\times V$,
we denote $S_{\mathrm{sym}}:=\{(x,y),(y,x)\mid(x,y)\in S\}$ and by
$S_{\mathrm{sym}}^c$ its complement in $V\times V$. The graph
$G\setminus S:=\bigl(E|_{S_{\mathrm{sym}}^c},V\bigr)$ is obtained by
deleting the edges in $S_{\mathrm{sym}}$ from~$G$.

\proof[Proof of Theorem~\ref{t:main0}]
Suppose that $G=(E,V)$ is a locally finite tree with bounded weights.
In view of \eqref{e:relation}, it is enough to consider a discrete
Schr\"odinger operator~$\Hc_G$ of the form $\Ac_G+\Vc$
for some potential $\Vc:V\to\R$.
Suppose that $\Hc_G$ has \underline{finite} and
positive deficiency index $\eta(\Hc_G)>0$.
Given a subgraph~$G'=(E', V')$ of~$G$,
we denote by $\Hc_{G'}$ the Schr\"odinger operator given by
$\Ac_{G'}+\Vc|_{V'}$.

We construct inductively a sequence $(v_k)_{k\in\N}$ of points of~$V$,
so that $v_k\sim v_{k+1}$ for all $k\in\N$.
Along the way we also define a sequence of subgraphs~$(G_k)_{k\in\N}$ of~$G$,
such that $G_k:=(E_k,V_k)$ is a tree, satisfying $v_k\in V_k$
and $\eta(\Hc_{G_k})\ge\eta(\Hc_{G_{k+1}})>0$.
Start with $G_0:=G$ and some $v_0\in V$.
For each $k\in\N$ we first remove the edges connected to~$v_k$
and obtain $G'_k:=G_k\setminus(\{v_k\}\times\Nr_{G_k}(v_k))$.
Using Lemma \ref{l:surgery} and the fact that $G_k$ is a tree, we find
\begin{equation*}
  0<\eta(\Hc_{G_k})=\eta(\Hc_{G'_k})
    =\sum\nolimits_{w\in\Nr_{G_k}(v_k)}\eta(\Hc_{\cc_{G'_k}(w)})\text.
\end{equation*}
Therefore there exists $w\in\Nr_{G_k}(v_k)$ with $\eta(\Hc_{\cc_{G'_k}(w)})>0$.
Set $v_{k+1}:=w$ and $G_{k+1}:=\cc_{G'_k}(w)$.
As announced the graph $G_{k+1}$ is a tree.

Since $k\mapsto\eta(\Hc_{G_k})$ is decreasing, positive, and has integer values,
there is $K\in\N$ so that $\eta(\Hc_{G_k})$
is constant for all $k\ge K$. Now consider
$\tilde G_k:=G_k[V_k\setminus V_{k+1}]$.

\begin{eqnarray*}%[h]
% \shorthandoff{"}	% xy braucht das Zeichen ", nur bei deutscher Umgebung
%  \begin{center}
%    \leavevmode
  \begin{xy}
    \xygraph {%
      %*{\bullet}-\ar@*{\text A}
      *{}-@{--}[r]*{\bullet}="A1"
	"A1"-@//[uul]*{\bullet}-@{.}[u]
	"A1"-@//[uu]*{\bullet}-@{.}[u]
	"A1"-@//[uur]*{\bullet}-@{.}[u]
	"A1"[ul]*{\tilde G_1}
	"A1"!{"A1"+<0cm,-3mm>*{v_1}}
      "A1"-@{--}[rrr]*{\bullet}="A1"
	-@// [u]*{\bullet}="B1"
	"B1"-@//[r]*{\bullet}%-@{-}[u]*{\bullet}-@{.}[u]
	"B1"-@//[u]*{\bullet}-@{.}[u]
	"B1"-@//[ul]*{\bullet}-@{.}[u]
	"A1"[ul]*{\tilde G_2}
	"A1"!{"A1"+<0cm,-3mm>*{v_2}}
      "A1"-@{--}[rrr]*{\bullet}="A1"
	-@// [u]*{\bullet}="B1"
	"B1"-@//[ru]*{\bullet}-@{.}[u]
	"B1"-@//[r]*{\bullet}
	"B1"-@//[u]*{\bullet}="B2"-@{.}[u]
	"B2"-@//[l]*{\bullet}-@{.}[u]
	"A1"[ul]*{\tilde G_3}
	"A1"!{"A1"+<0cm,-3mm>*{v_3}}
      "A1"-@{--}[rrr]*{\bullet}="A1"
	-@//[ul]*{\bullet}-@//[u]*{\bullet}-@{.}[u]
	"A1"-@//[u]*{\bullet}="B1"
	"B1"-@//[ru]*{\bullet}-@{.}[u]
	"B1"-@//[u]*{\bullet}-@{.}[u]
	"A1"[ur]*{\tilde G_4}
	"A1"!{"A1"+<0cm,-3mm>*{v_4}}
      "A1"-@{--}[r]*{}
    }%
  \end{xy}
%  \end{center}%
% \shorthandon{"}% Nur bei deutscher Umgebung
\end{eqnarray*}%
\emph{The dashed line at the bottom shows the constructed path~$(v_k)_{k\in\N}$.
  The graphs~$\tilde G_k$ can extend infinitely, as indicated with the dots.
  For all $k\in\N$ the graph~$G_k$ is the union of all~$\tilde G_{k'}$
  with $k'\ge k$ plus the dashed bottom line starting at~$v_k$.
}%
\label{f:labelsintree}%

Again by Lemma~\ref{l:surgery}, we infer
\begin{equation*}
\eta(\Hc_{\tilde G_k})=\eta(\Hc_{G_k})-\eta(\Hc_{G_{k+1}})=0,
\text{ for $k\ge K$.}
\end{equation*}
By one more application of Lemma~\ref{l:surgery}, we obtain
\begin{equation*}
0<\eta(\Hc_{G_K})
  =\sum_{k=K}^\infty\eta(\Hc_{\tilde G_k})
  =0\text.
\end{equation*}
This is a contradiction.\qed

\section{Recursive graphs}\label{s:recu}
In this section we generalize the previous approach to graphs
which satisfy a certain recursive property.
We shall use the notation related to subgraphs, which were introduced
in section \ref{s:possible}. To simplify notation, given a potential
$\Vc:V\to\R$ and the Schr\"odinger operator $\Hc_G:=\Ac_G+\Vc$
acting on $G=(E,V)$, we shall write $\eta_\Hc(G):=\eta(\Hc_G)$.
As above, if $G'$ is a subgraph of~$G$ obtained by removing edges,
we denote by $\Hc_{G'}$ the Schr\"odinger operator $\Ac_G+\Vc$.

\begin{definition}\label{d:r}
Let $G=(E,V)$ be a locally finite graph.
Given $M>0$ and $\Vc:V\to \R$, we say that $G$ has the
property~$\Rc(M,\Vc)$, if
\begin{itemize}
  \item
    either $\eta_\Hc(G)=0$ or
  \item we can find a partition $\{B,U_n,W_n\mid n\in\N\}$ of~$V$
      such that
    \begin{enumerate}[\rm(P1)]
      \item $\eta_\Hc({G[B]})=0$,
      \item $\{\tilde U_n,\tilde W_n\mid n\in\N\}$ is pairwise disjoint,
	where $\tilde U_n:=B\cap\Nr_G(U_n)$ and $\tilde
        W_n:=B\cap\Nr_G(W_n)$.
      \item For $m,n\in\N$, $E(U_n,W_m)=0$ and
	$E(U_n,U_m)=0$, $E(W_n,W_m)=0$ if $m\ne n$,
      \item $\forall x\in B\colon\lvert\Nr_G(x)\cap U_n\rvert\le M$ and
	$\forall x\in U_n\colon\lvert\Nr_G(x)\cap B\rvert\le M$ for all $n\in\N$,
      \item $\forall x\in\tilde W_n\colon
	 \lvert\Nr_G(x)\cap(B\setminus\tilde W_n)\rvert\le M$ and
	$\forall x\in\tilde B\setminus\tilde W_n\colon
	 \lvert\Nr_G(x)\cap B\rvert\le M$ for all $n\in\N$,
      \item $G[U_n]$ and $G[W_n\cup\tilde W_n]$
	have the property~$\Rc(M,\Vc|_{U_n})$
	and $\Rc\big(M,\Vc|_{W_n\cup\tilde W_n}\big)$, respectively.
    \end{enumerate}
\end{itemize}
\end{definition}
We explain in words, what the sets~$B$, $U_n$, and $W_n$ are.
The set of vertices $B\subseteq V$ stands for the \emph{base}
of the graph~$G$.  We recall that, by definition, $G[B]$ is the restriction of
the graph~$G$ to~$B$, see the Introduction.  In the case of trees,
for the $k$-th level of recursion we use $B=\{v_k\}$.
We allow more complicated situations here:
for instance, $\eta_\Hc(B)=0$ if $B$ has bounded degree and weights,
see also Proposition \ref{p:essSA}.
The graphs $G[U_n]$ and $G[W_n\cup\tilde W_n]$
correspond to subgraphs that we want to cut out
and to study in the next recursive step. The condition (P4) ensures that each
element of the subgraph~$G[U_n]$ is linked by at most~$M$ edges to the base.
On the other hand, the graph $G[W_n]$ could be linked to the base by a
large number of edges, like in the previous case for trees.
In this situation, we shall not consider $G[W_n]$
in the next recursive step but $G[W_n\cup\tilde W_n]$,
which contains a part of the base. Notice that
$B\setminus\bigcup\nolimits_{n\in\N}\tilde W_n$ is empty in the
previous setting of a tree.
Note that condition (P5) makes sure that each element of the subgraph
$G[W_n\cup\tilde W_n]$ is linked to the remaining part of the base
with at most $M$~edges.
Condition (P2) ensures that the subgraphs~$G[U_n]$
and $G[W_n\cup\tilde W_n]$ are not too close to each other.
Condition (P3) tells that there are no edges between the $U_n$ and $W_n$.
This condition can be relaxed with Lemma~\ref{l:surgery},
by asking that each vertices is linked with at most $M$ other ones.

Definition~\ref{d:r} is motivated by the fact that the recursive process
splits the deficiency indices in a conservative way.
\begin{lemma}\label{l:index}
  Suppose that $G$ is a locally finite graph with bounded weights
  satisfying \textrm{(P2) -- (P5)}.
  Then, using the notation of Definition~\ref{d:r},
  \begin{equation}\label{e:index1}
    \eta_{\Hc}(G)
      =\eta_\Hc\left({G\bigl[B\setminus\bigcup\nolimits_{n\in\N}\tilde
	W_n\bigr]}\right)+\sum_{n\in\N}\eta_\Hc\bigl(G[U_n]\bigr)
	  +\eta_\Hc\bigl(G[W_n\cup\tilde W_n]\bigr).
  \end{equation}
  Moreover, if $G$ obeys \textrm{(P1)},
  \begin{equation}\label{e:index2}
    \eta_\Hc\left({G\bigl[B\setminus\bigcup\nolimits_{n\in\N}\tilde
      W_n\bigr]}\right)=0.
  \end{equation}
\end{lemma}
\proof
Equation~\eqref{e:index1} is a direct consequence of Lemma~\ref{l:surgery}.
By the same argument
\begin{equation*}
0=\eta_\Hc(B)=
\eta_\Hc\left({G\bigl[B\setminus\bigcup\nolimits_{n\in\N}\tilde
    W_n\bigr]}\right)+\sum_{n\in\N}\eta_\Hc\left({G[\tilde W_n]}\right).
\end{equation*}
Equation~\eqref{e:index2} follows, since deficiency indices are nonnegative.
\qed

Finally, we prove:
\begin{theorem}\label{t:main}
Suppose that $G$ is a locally finite graph with bounded weights
satisfying property $\Rc(M,\Vc)$, for a certain potential~$\Vc$.
Then $\eta(\Ac_G+\Vc)= \eta_\Hc(G)\in \{0, \infty\}$.
\end{theorem}
\proof
Let $G$ be a graph fulfilling all assumptions and having
\underline{finite} and positive deficiency index.
As in the case of trees we construct a sequence of nested subgraphs
$(G_k)_{k\in\N}$ of~$G$ such that for all $k\in\N$
\begin{itemize}
\item $\eta_\Hc(G_k)\ge\eta_\Hc(G_{k+1})>0$,
\item $G_k$ satisfies property~$\Rc(M)$.
\end{itemize}
We set $G_0:=G$ and construct $G_{k+1}$ inductively from $G_k$.
We use now Lemma~\ref{l:index}.  Taking advantage of \eqref{e:index2}
in \eqref{e:index1}, there is a subgraph of $G_k$, among the family
$\{G_k[U_n(k)]$, $G_k[W_n(k)\cup\tilde W_n(k)]\mid k\in\N\}$
with positive deficiency index.
We call it $G_{k+1}$.
By \eqref{e:index1} we have $\eta_\Hc(G_k)\ge\eta_\Hc(G_{k+1})$.
Thanks to {\rm(P5)}, $G_{k+1}$ satisfies also property~$\Rc(M)$.

As in Theorem~\ref{t:main0} we conclude that
there is $K\in\N$ so that $\eta_\Hc({G_k})$
is constant for all $k\ge K$.
Now consider $\tilde G_k:=G_k[V_k\setminus V_{k+1}]$.
By Lemma~\ref{l:surgery}, we infer
$\eta_\Hc({\tilde G_k})=\eta_\Hc({G_k})-\eta_\Hc({G_{k+1}})=0$, for $k\ge K$.
By construction there are at most $M$
connections per vertex between $G_k$ and $G_{k+1}$.
By a last application of Lemma~\ref{l:surgery}, we obtain
$0<\eta_\Hc({G_K})=\sum_{k=K}^\infty\eta_\Hc({\tilde G_k})=0$.
This is the desired contradiction.
\qed

We finish by mentioning a possible generalization.

\begin{remark} In the previous result, we do not suppose more than
having bounded weights. The main examples we have in mind are simple graphs.
However if one considers weighted graphs such that
$\inf\bigl(E(V\times V)\setminus\{0\}\bigr)=0$,
using \eqref{e:surgery},
one can relax the hypothesis on the uniformity in~$M$,
which is implemented in Definition~\ref{d:r}.
\end{remark}
\appendix
\renewcommand{\theequation}{\thesection .\arabic{equation}}

\section{Stability of the deficiency indices of a symmetric
  operator}\label{s:stabind}
\setcounter{equation}{0}

Given a closed and densely defined symmetric operator~$S$, one has the
obvious inclusion $\Dc(S)\subset\Dc(S^*)$.
In fact, given $z\in\C\setminus\R$, one gets the topological direct sum
\begin{eqnarray}\label{e:dom}
  \Dc(S^*)=\Dc(S)\oplus\ker(S^*+z)\oplus\ker(S^*-\overline{z}).
\end{eqnarray}
One also knows that $z\mapsto \dim\big(\ker(S^*-z)\big)$ is constant on the two
connected components of $\C\setminus \R$. Note also that $\dim
\left(\Dc(S^*)/\Dc(S)\right)= \eta_-(S)+\eta_+(S)$.
We refer to \cite[Section X.1]{RS} for an introduction to the subject.

For the convenience of the reader and as we were not able to
locate a proof in the literature, we recall the following useful and
well-known fact. It is essentially due to Kato and Rellich.
\begin{proposition}\label{p:stab}
  Given two closed and densely defined symmetric operators $S,T$
  acting on a complex Hilbert space and such that $\Dc(S)\subset\Dc(T)$.
  Suppose there are $a\in[0,1)$ and $b\ge0$ such that
  \begin{eqnarray}\label{e:KR}
    \norm{Tf}\le a\norm{Sf}+b\norm{f},\text{ for all }f\in\Dc(S).
  \end{eqnarray}
  Then, the closure of $(S+T)|_{\Dc(S)}$ is a symmetric operator that
  we denote by $S+T$.
  Moreover, one obtains that $\Dc(S)=\Dc(S+T)$
  and that $\eta_\pm(S)=\eta_\pm(S+T)$.
  In particular, $S+T$ is self-adjoint if and only if~$S$ is.
\end{proposition}
Note that if~$S$ is self-adjoint, i.e., $\eta_\pm(S)=0$, the above
result is the standard Kato-Rellich theorem, e.g., \cite[Theorem
X.12]{RS}. In the proofs of this article, we use this result in the
case $a=0$ and $\eta_-(S)= \eta_+(S)$.
In this setting, one can avoid this general result and
repeat a shorter argumentation, coming from \cite{Gol}. We explain
this alternative approach at the end of the proof of Lemma
\ref{l:surgery}. Finally, we point out that all the results about
deficiency indices of this article are stable under the above class of
perturbation, i.e., \eqref{e:KR} with $a\in [0,1)$.

\proof
Let $\theta\in[-1,1]$. Note that $W_\theta|_{\Dc(S)}:=\left(
S+\theta T\right)|_{\Dc(S)}$ is symmetric and closable. Its closure
is denoted by $W_\theta$. Using \eqref{e:KR}, one sees that the graph
norms of~$S$ and of $W_\theta$ are equivalent on $\Dc(S)$. Then,
we infer that $W_\theta$ is closed, symmetric and with domain
$\Dc(W_\theta)=\Dc(S)$. In particular, $\Dc(S+T)=\Dc(S)$.

We concentrate on the deficiency indices.
It is enough to consider the case $a\in(0,1)$ and $b>0$.
Notice first that, for $f\in\Dc(S)$ and $\varepsilon>0$, one obtains
$\norm{Tf}^2\le a^2(1+\varepsilon)\norm{Sf}^2+b^2(1+1/\varepsilon)\norm{f}^2$
for all $\varepsilon>0$.
Then, since~$S$ is symmetric, we derive that
\begin{eqnarray}\label{e:alpha}
\norm{Tf}^2\leq\alpha^2\norm{(S\pm\rmi\gamma)f}^2, \mbox{ for all
}f\in\Dc(S)
\end{eqnarray}
and where $\alpha^2=(1+\varepsilon)a^2$
and $\gamma=\sqrt{b^2/(\varepsilon a^2)}$.
Taking $\varepsilon$ small enough,
we reduce to the case $\alpha\in(0,1)$ and $\gamma\geq1$.
Take now
\begin{eqnarray}\label{e:theta}
\theta_1, \theta_2\in [-1,1], \mbox{ so that }
|\theta_1-\theta_2|<\frac{(1-\alpha)}{\alpha}.
\end{eqnarray}
We now prove:
\begin{eqnarray}\label{e:ker}
 \ker((W_{\theta_1})^*\pm \rmi \gamma)\cap
\big(\ker((W_{\theta_2})^*\pm\rmi \gamma)\big)^\perp=
\ker((W_{\theta_1})^*\pm \rmi \gamma)\cap
\ran(W_{\theta_2}\mp\rmi \gamma)= \{0\}.
\end{eqnarray}
Given~$H$ a closed symmetric and densely defined operator,
by considering $\im \langle y, (H\pm\rmi \gamma)
y\rangle$, one sees that $\norm{(H \pm\rmi \gamma) y}\geq \gamma \norm{y}$
for all $y\in \Dc(H)$ and that the range of $(H \pm\rmi \gamma)$ is
closed. Hence, the first equality of \eqref{e:ker} holds true.

Take $x\in\Dc\big((W_{\theta_1})^*\big)\setminus\{0\}$
and in the intersection in the l.h.s.\ of \eqref{e:ker}.
We finish the proof for the minus sign. The other case is done analogous.
We infer that there is $z\in \Dc(S)\setminus\{0\}$,
such that $(W_{\theta_2}+\rmi \gamma)z=x$. Then,
\begin{align}\label{e:newequa}
0&= \langle ( (W_{\theta_1})^* -\rmi \gamma) x, z\rangle= \langle x,
(W_{\theta_1}+\rmi\gamma )z\rangle = \norm{x}^2 + (\theta_2-\theta_1)
\langle x, T z\rangle.
\end{align}
Now, with \eqref{e:alpha}, we infer $(1-\alpha)\norm{Tz}\leq
\alpha\norm{(W_{\theta_2}+ \rmi \gamma) z}$. Using the latter with
\eqref{e:theta} and \eqref{e:newequa}, we derive:
\[ \norm{x}\leq |\theta_1-\theta_2| \cdot \norm{T z}<
\norm{ (W_{\theta_2}+\rmi \gamma)z}= \norm{x},\]
which is a contradiction. This proves \eqref{e:ker} and therefore
$\dim\ker({W_{\theta_2}}^*\pm\rmi\gamma) \ge
\dim\ker((W_{\theta_1})^*\pm\rmi\gamma)$, under the hypothesis \eqref{e:theta}.
One deduces easily that $\dim\ker((W_\theta)^*\pm\rmi\gamma) =
\dim\ker(S^*\pm\rmi\gamma)$, for all $\theta\in [-1,1]$.
To conclude, we recall that, given a symmetric operator~$H$,
one has that $z\mapsto \dim\big(\ker(H^*-z)\big)$
is constant on the two connected components of $\C\setminus\R$.
\qed

We now give a direct application to Jacobi matrices, which act on
$\ell^2(\N)$. Given~$A$, the closure of a three-diagonal symmetric
Jacobi matrix with $a_n\in \R$ on the diagonal and $b_n>0$ on the
upper diagonal, it is well known, e.g., \cite[Page 504]{Be},
that if $\sum_{n\in\N}1/b_n=\infty$ and with no condition on
the sequence $(a_n)_{n}$, then $A^*=A$.
We give a generalisation in Proposition~\ref{p:essSA} \textit{(2)}.
With again no condition on the diagonal elements, we prove now:

\begin{proposition}\label{p:jacobi}
  Let~$A$ be the closure of a $(2N+1)$-diagonal (complex-)symmetric matrix
  acting by $Af(n)=\sum_{k\in\N}a_{k,n}f(k)$ for $f:\N\to\C$
  with compact support and where $a_{k,n}\in\C$, for $k,n\in\N$. If
  \begin{eqnarray*}
    \liminf_{n\to\infty}c_n<\infty,
    \text{ where }c_n:=\max_{0\leq l<k\leq K}|a_{n-1-l, n-l+k}|,
  \end{eqnarray*}
  for $n\geq K$, then $A=A^*$.
\end{proposition}
\proof
Let $(c_{u_n})_{n\in\N}$ be a \emph{bounded} subsequence of $(c_n)_{n\in\N}$
and set $B_n:=\bone_{[u_n,u_{n+1}-1]}\,A\,\bone_{[u_n,u_{n+1}-1]}$.
\begin{eqnarray*}
  A=\left(\begin{array}{l}
  \xymatrix @!=0.25pc { \ar@{.}[dr]&
  \\
  & \ar@{-}[rr]\ar@{-}[dd]& & \ar@{-}[dd]
  \\
   & &{\mathbf {B_{n-1}}} & & \ar@{.}[d]\ar@{.}[dr]
  \\
  &\ar@{-}[rr] & & a_{u_{n}-1, u_{n}-1} & \ar@{.}[r]&
  \\&&\ar@{.}[dr]& \ar@{.}[d] \ar@{.}[l]& a_{u_n, u_n}\ar@{-}[rr]\ar@{-}[dd]& & \ar@{-}[dd]
  \\&&& & &\mathbf {B_{n}} &
  \\
  && &&\ar@{-}[rr] & &
  \\
  &&&&&&& \ar@{.}[ul]
  }
  \end{array}\right)
\end{eqnarray*}
Set~$B$ be the closure of $\oplus_n B_n$. Note that the
deficiency indices of~$B$ are $(0,0)$, since $B_n$ are finite
dimensional matrices. Then, remembering that $\sup_{n\in\N}|c_{u_n}|<\infty$,
we see that $(B-A)|_{\Cc_c(\N)}$ extends to a bounded operator.
 Therefore, Proposition~\ref{p:stab} entails that~$A$ is self-adjoint.
\qed

\end{document}